%&amstex
\input amsppt.sty
\magnification=\magstep1
\hsize=30truecc
\baselineskip=16truept
\vsize=22.2truecm
%\NoBlackBoxes
\nologo
\pageno=1
\topmatter
\TagsOnRight
\def\N{\Bbb N}
\def\Z{\Bbb Z}

\def\l{\left}
\def\r{\right}
\def\b{\bigg}

\def\({\b(}
\def\[{\b[}
\def\){\b)}
\def\]{\b]}

\def\t{\text}
\def\f{\frac}
\def\mo{\roman{mod}}

\def\bi{\binom}
\def\eq{\equiv}

\def\ls{\leqslant}

\def\al{\alpha}

\def\Proof{\noindent{\it Proof}}
\def\Remark{\noindent{\it Remark}}

\def\Ack{\noindent {\bf Acknowledgment}}
\hbox {Final version for publication in J. Number Theory.}
\medskip
\title On Euler numbers modulo powers of two\endtitle
\author Zhi-Wei Sun\endauthor
\affil Department of Mathematics (and Institute of Mathematical Science)
\\Nanjing University, Nanjing 210093, The People's Republic of China
\\  \tt{zwsun\@nju.edu.cn}
\\ {\tt http://pweb.nju.edu.cn/zwsun}
\endaffil
\abstract To determine Euler numbers modulo powers of two seems to be a
difficult task. In this paper we achieve this and apply the explicit congruence
to give a new proof of a classical result due to M. A. Stern.
\endabstract
\thanks 2000 {\it Mathematics Subject Classification}. Primary 11B68;
Secondary 11A07, 11S05.
\newline\indent
The author is supported by the National Science Fund for Distinguished Young Scholars (No. 10425103)
and the Key Program of NSF (No. 10331020) in P. R. China.
\endthanks
\endtopmatter
\document

\heading {1. Introduction}\endheading
Euler numbers $E_0,E_1,E_2,\ldots$ are integers given by
$$E_0=1\ \ \t{and}\ \ E_n=-\sum^{n-1}\Sb k=0\\2\mid n-k\endSb\bi nkE_k\ \t{for}\ n\in\Z^+=\{1,2,3,\ldots\}.$$
It is well known that $E_{2n+1}=0$ for each $n\in\N=\{0,1,2,\ldots\}$ and
$$\sec x=\sum_{n=0}^{\infty}(-1)^nE_{2n}\f{x^{2n}}{(2n)!}\ \ \l(|x|<\f{\pi}2\r).$$
The Euler polynomial $E_n(x)$ of
degree $n$ is defined by
$$E_n(x)=\sum_{k=0}^n\bi nk\f {E_k}{2^k}\l(x-\f12\r)^{n-k}.$$
Clearly $E_n=2^nE_n(1/2)$. That $E_n(x)+E_n(x+1)=2x^n$ is a well-known fact.
[Su3] contains some symmetric identities for Euler polynomials.

Euler numbers modulo an odd integer are trivial. In fact, for any
$k\in\N$ and $q\in\Z^+$ we have
$$2^kE_k\l(q+\f12\r)=2^k\sum_{l=0}^k\bi kl\f{E_l}{2^l}q^{k-l}\eq E_k=2^kE_k\l(\f12\r)\ (\mo\ q)$$
and
$$\aligned &E_k\l(\f12\r)-(-1)^qE_k\l(q+\f12\r)
\\=&\sum_{j=0}^{q-1}
\((-1)^jE_k\l(j+\f12\r)-(-1)^{j+1}E_k\l(j+1+\f12\r)\)
\\=&2\sum_{j=0}^{q-1}(-1)^j\l(j+\f12\r)^k,
\endaligned$$
therefore
$$E_k\eq\sum_{j=0}^{q-1}(-1)^j(2j+1)^k\ (\mo\ q)\quad \ \t{providing}\ 2\nmid q.\tag1.1$$

 It is natural to determine
Euler numbers modulo powers of two.
However, this is a difficult task since $1/2$ is not a $2$-adic integer.
As far as the author knows, no one else has achieved this before.

 In this paper we determine Euler numbers modulo powers of two in the following explicit way.

\proclaim{Theorem 1.1} Let $n\in\Z^+$. If $k\in\N$ is even, then
$$\f{3^{k+1}+1}4E_k\eq\f{3^k}2\sum_{j=0}^{2^n-1}
(-1)^{j-1}(2j+1)^k\l\lfloor\f{3j+1}{2^n}\r\rfloor\ (\mo\ 2^n)\tag1.2$$
where $\lfloor\al\rfloor$ denotes the greatest integer not exceeding a real number $\al$,
moreover for any positive odd integer $m$ we have the congruence
$$\aligned&\f{m^{k+1}-(-1)^{(m-1)/2}}4E_k
\\\eq&\f{m^k}2\sum_{j=0}^{2^n-1}(-1)^{j-1}(2j+1)^k\l\lfloor\f{jm+(m-1)/2}{2^n}\r\rfloor
\ (\mo\ 2^n).
\endaligned\tag1.3$$
\endproclaim

\Remark\ 1.1. When $k\in\N$ is even, $3^{k+1}+1\eq3+1=4\ (\mo\ 8)$ and so $(3^{k+1}+1)/4$ is odd.
\medskip

Theorem 1.1 implies the following nice result.

\proclaim{Theorem 1.2} Let $k,l\in\N$ be even. If $2^n\|(k-l)\ (\t{i.e.}, 2^n|(k-l)$
but $2^{n+1}\nmid (k-l))$
where $n\in\Z^+$, then $2^n\|(E_k-E_l)$. In other words, for any $n\in\Z^+$ we have
$$E_k\eq E_l\ (\mo\ 2^n)\iff k\eq l\ (\mo\ 2^n).\tag1.4$$
\endproclaim

\Remark\ 1.2. Theorem 1.2 does not
tell us how to determine $E_k$ mod $2^n$ with $0\ls k<2^n$.
In 1875 Stern [S] gave a brief sketch of a proof of Theorem 1.2,
then Frobenius amplified Stern's sketch in 1910. In 1979 Ernvall [E] said that he could not understand
Frobenius' proof and provided his own proof
involving umbral calculus.
Recently an induction proof of Theorem 1.2 was given by Wagstaff [W].
\medskip

In the next section we will provide some lemmas.
Theorems 1.1 and 1.2 will be proved in Section 3.

Now we introduce some notations throughout this paper.
For $a,b\in\Z$ by $(a,b)$ we mean the
greatest common divisor of $a$ and $b$.
For an integer $q>1$, we use $\Z_q$ to denote the ring
of rational $q$-adic integers (see [M] for an introduction to $q$-adic numbers);
those numbers in $\Z_q$ are simply called {\it $q$-integers}
and they have the form $a/b$ with $a\in\Z$, $b\in\Z^+$ and $(a,b)=(b,q)=1$.
For $\al,\beta\in\Z_q$, by $\al\eq\beta\
(\mo\ q)$ we mean that $\al-\beta=q\gamma$ for some
$\gamma\in\Z_q$. Two polynomials in $\Z_q[x]$ are said to be congruent modulo $q$
if all the corresponding coefficients in the two polynomials are congruent mod $q$.

Bernoulli numbers $B_0,B_1,B_2,\ldots$ given by
$B_0=1$ and the recursion
$$\sum_{k=0}^n\bi{n+1}kB_k=0\ \ (n=1,2,3,\ldots)$$
are closely related to Euler numbers.

Let $k\in\Z^+$ be even, and let $p$ be an odd prime with $p-1\nmid k$.
In 1851 E. Kummer showed that $B_k/k\in\Z_p$, and $B_k/k\eq B_l/l\ (\mo\ p^n)$
(where $n\in\Z^+$)
for any $l\in\Z^+$ with $k\eq l\ (\mo\ \varphi(p^{n}))$,
where $\varphi$ is Euler's totient function.
In contrast with Theorem 1.2, the converse of Kummer's congruences is not true, e.g.,
$$\f{B_{16}}{16}\eq \f{B_4}4\ (\mo\ 13^2)\ \ \t{but}\ 16\not\eq4\ (\mo\ \varphi(13^2)).$$

Suppose that $p^n\|k$ where $n\in\N$. Then $p^n$ divides the numerator of $B_k$ since
$p$ does not divide the denominator
of $B_k$ by the von Staudt--Clausen theorem (cf. [IR, p.\,233]). In [T] this trivial observation
was attributed to J. C. Adams. Recently R. Thangadurai [T] conjectured that if $n>0$ then
$p^{n+2}$ does not divide the numerator of $B_k$ (i.e., $B_k/k\not\in p^2\Z_p$).

\heading{2. Several Lemmas}\endheading

For each $n\in\N$ the Bernoulli polynomial $B_n(x)$ of degree $n$ is given by
$$B_n(x)=\sum_{k=0}^n\bi nkB_kx^{n-k}.$$
A useful multiplication formula of Raabe asserts that
$$m^{n-1}\sum_{r=0}^{m-1}B_n\l(\f{x+r}m\r)=B_n(x)\ \ \t{for any}\ m\in\Z^+.$$
Euler polynomials are related to Bernoulli polynomials in the following manner:
$$\aligned \f{n+1}2E_n(x)=&B_{n+1}(x)-2^{n+1}B_{n+1}\l(\f x2\r)
\\=&2^{n+1}B_{n+1}\l(\f{x+1}2\r)-B_{n+1}(x).\endaligned\tag2.1$$

\proclaim{Lemma 2.1 {\rm ([Su2, Cor. 1.3])}} Let $a\in\Z$ and $k,m\in\Z^+$.
Let $q>1$ be an integer relatively prime to $m$. Then
$$\aligned&\f1k\(m^kB_k\l(\f{x+a}m\r)-B_k(x)\)
\\\eq&\sum_{j=0}^{q-1}\(\l\lfloor\f{a+jm}q\r\rfloor+\f{1-m}2\)(x+a+jm)^{k-1}\ (\mo\ q).
\endaligned\tag2.2$$
\endproclaim

\Remark\ 2.1. If $q>1$ is an integer relatively prime to $m\in\Z$,
then $(1-m)/2\in\Z_q$ because $q$ or $m$ is odd.

\proclaim{Lemma 2.2} Let $a\in\Z$, $k\in\N$ and $m\in\Z^+$.
Let $q\in\Z^+$, $2\mid q$ and $(m,q)=1$. Then
$$\aligned&\f{m^{k+1}}2E_k\l(\f{x+a}m\r)-\f{(-1)^a}2E_k(x)
\\\eq&\sum_{j=0}^{q-1}(-1)^{j-1}\(\l\lfloor\f{a+jm}q\r\rfloor+\f{1-m}2\)(x+a+jm)^{k}\ (\mo\ q).
\endaligned\tag2.3$$
\endproclaim
\Proof. Let us first handle the case $4\mid q$.
Denote by $\{a\}_2$ the least nonnegative residue of $a$ modulo $2$ and set
$\bar x=(x+\{a\}_2)/2$. In view of (2.1) we have
$$\align&\f{m^{k+1}}2E_k\l(\f{x+a}m\r)-\f{(-1)^a}2E_k(x)
\\=&\f{m^{k+1}}{k+1}\(B_{k+1}\l(\f{x+a}m\r)-2^{k+1}B_{k+1}\l(\f{x+a}{2m}\r)\)
\\&-\f{1}{k+1}\l(B_{k+1}(x)-2^{k+1}B_{k+1}(\bar x)\r)
\\=&\f{1}{k+1}\(m^{k+1}B_{k+1}\l(\f{x+a}m\r)-B_{k+1}(x)\)
\\&-\f{2^{k+1}}{k+1}\(m^{k+1}B_{k+1}\l(\f{\bar x+\lfloor a/2\rfloor}m\r)-B_{k+1}(\bar x)\).
\endalign$$

Let $P(t)$ denote the polynomial
$$\align&\f1{k+1}\(m^{k+1}B_{k+1}\l(\f{t+\lfloor a/2\rfloor}m\r)-B_{k+1}(t)\)
\\-&\sum_{j=0}^{q/2-1}\(\l\lfloor\f{\lfloor a/2\rfloor+jm}{q/2}\r\rfloor+\f{1-m}2\)
\l(t+\l\lfloor\f a2\r\rfloor+jm\r)^k.
\endalign$$
Clearly $\deg P(t)\ls k$. Recall that $4\mid q$. By Lemma 2.1,
we can write
$$P(t)=\sum_{i=0}^k\f q2c_it^i\quad\t{where}\ c_i\in\Z_{q/2}=\Z_q.$$
Thus
$$\f{2^{k+1}P(\bar x)}q=\sum_{i=0}^kc_i2^k\l(\f{x+\{a\}_2}2\r)^i\in\Z_q[x].$$

In light of the above,
$$\align&\f{m^{k+1}}2E_k\l(\f{x+a}m\r)-\f{(-1)^a}2E_k(x)
\\\eq&\f1{k+1}\(m^{k+1}B_{k+1}\l(\f{x+a}m\r)-B_{k+1}(x)\)
\\&-2^{k+1}\sum_{j=0}^{q/2-1}\(\l\lfloor\f{a/2+jm}{q/2}\r\rfloor+\f{1-m}2\)
\l(\bar x+\l\lfloor\f a2\r\rfloor+jm\r)^k
\\\eq&\sum_{j=0}^{q-1}\(\l\lfloor\f{a+jm}q\r\rfloor+\f{1-m}2\)(x+a+jm)^k
\\&-2\sum^{q-1}\Sb i=0\\2\mid i\endSb\(\l\lfloor\f{a+im}q\r\rfloor+\f{1-m}2\)(x+a+im)^k
\\\eq&\sum_{j=0}^{q-1}(-1)^{j-1}\(\l\lfloor\f{a+jm}q\r\rfloor+\f{1-m}2\)(x+a+jm)^k\ (\mo\ q).
\endalign$$
Therefore (2.3) holds.

Now we consider the remaining case $2\|q$. As $4\mid 2q$ and $(m,2q)=1$, by the above we have
$$\align&\f{m^{k+1}}2E_k\l(\f{x+a}m\r)-\f{(-1)^a}2E_k(x)
\\\eq&\sum_{j=0}^{2q-1}(-1)^{j-1}\(\l\lfloor\f{a+jm}{2q}\r\rfloor+\f{1-m}2\)(x+a+jm)^k\ (\mo\ 2q)
\\\eq&\sum_{j=0}^{q-1}(-1)^{j-1}\(\l\lfloor\f{a+jm}{2q}\r\rfloor+\f{1-m}2\)(x+a+jm)^k
\\&+\sum_{j=0}^{q-1}(-1)^{j+q-1}\(\l\lfloor\f{a+(j+q)m}{2q}\r\rfloor+\f{1-m}2\)(x+a+jm)^k\ (\mo\ q)
\\\eq&\sum_{j=0}^{q-1}(-1)^{j-1}\l(a_j+\f{1-m}2\r)(x+a+jm)^k\ \ (\mo\ q),
\endalign$$
where
$$\align a_j=&\l\lfloor\f{a+jm}{2q}\r\rfloor+\l\lfloor\f{a+jm+q}{2q}\r\rfloor
\\=&\l\lfloor\f{(a+jm)/q}{2}\r\rfloor+\l\lfloor\f{(a+jm)/q+1}{2}\r\rfloor
=\l\lfloor\f{a+jm}q\r\rfloor.
\endalign$$
This completes the proof. \qed

\proclaim{Lemma 2.3} Let $a\in\Z$, $m,q\in\Z^+$, $2\mid q$ and $(m,q)=1$. Then
$$\sum_{j=0}^{q-1}(-1)^{j-1}\(\l\lfloor\f{a+jm}q\r\rfloor+\f{1-m}2\)=\f{m-(-1)^a}2.\tag2.4$$
\endproclaim
\Proof. Observe that
$$\aligned&\sum_{j=0}^{q-1}(-1)^j\(\l\lfloor\f{a+jm}q\r\rfloor+\f{1-m}2\)
\\=&2\sum^{q-1}\Sb j=0\\2\mid j\endSb\(\l\lfloor\f{a+jm}q\r\rfloor+\f{1-m}2\)
-\sum_{j=0}^{q-1}\(\l\lfloor\f{a+jm}q\r\rfloor+\f{1-m}2\)
\\=&2\sum_{i=0}^{q/2-1}\(\l\lfloor\f{a/2+im}{q/2}\r\rfloor+\f{1-m}2\)
-\sum_{j=0}^{q-1}\(\l\lfloor\f{a+jm}q\r\rfloor+\f{1-m}2\)
\\=&2\l(\l\lfloor\f a2\r\rfloor+\f{1-m}2\r)-\l(\lfloor a\rfloor+\f{1-m}2\r)\ (\t{by [Su1, Prop. 2.1]})
\\=&2\l\lfloor\f a2\r\rfloor-a+\f{1-m}2=\f{(-1)^a-m}2.
\endaligned$$
This proves (2.4). \qed

\heading{3. Proofs of Theorems 1.1 and 1.2}\endheading

\noindent {\it Proof of Theorem 1.1}. Obviously (1.2) follows from (1.3) in the case $m=3$.

 Let $m\in\{1,3,5,\ldots\}$ and $a=(m-1)/2$.
Applying Lemma 2.2 with $q=2^n$, we find that the polynomial
$$\align f(x):=&\f{m^{k+1}}2E_k\l(\f{x+a}m\r)-\f{(-1)^a}2E_k(x)
\\&+\sum_{j=0}^{2^n-1}(-1)^{j}\(\l\lfloor\f{a+jm}{2^n}\r\rfloor+\f{1-m}2\)(x+a+jm)^k
\endalign$$
belongs to $2^n\Z_{2^n}[x]$. Observe that
the coefficient of $x^k$ in $f(x)$ is zero because
$$\f{m^{k+1}}2\l(\f1m\r)^k-\f{(-1)^a}2
+\sum_{j=0}^{2^n-1}(-1)^{j}\(\l\lfloor\f{a+jm}{2^n}\r\rfloor+\f{1-m}2\)=0$$
by Lemma 2.3. So $\deg f(x)\ls k-1$ and hence $2^{k-1}f(x/2)\in
2^n\Z_{2^n}[x]$. (Note that $f(x)=0$ if $k=0$.) In particular, $2^{k-1}f(1/2)\eq0\ (\mo\ 2^n)$.

Clearly
$$\align 2^{k-1}f\l(\f12\r)=&\f{m^{k+1}}2 2^{k-1}E_k\l(\f12\r)-\f{(-1)^a}2 2^{k-1}E_k\l(\f12\r)
\\&+2^{k-1}\sum_{j=0}^{2^n-1}(-1)^{j}\(\l\lfloor\f{a+jm}{2^n}\r\rfloor+\f{1-m}2\)\l(\f m2+jm\r)^k
\endalign$$
and thus
$$\align&2^{k-1}f\l(\f12\r)-\f{m^{k+1}-(-1)^{(m-1)/2}}4E_k
\\=&\f{m^k}2\sum_{j=0}^{2^n-1}(-1)^{j}\(\l\lfloor\f{jm+(m-1)/2}{2^n}\r\rfloor+\f{1-m}2\)(2j+1)^k.
\endalign$$
So it remains to show
$$\sum_{j=0}^{2^n-1}(-1)^j(2j+1)^k\eq0\ (\mo\ 2^{n+1})\quad\t{providing}\ 2\mid k.$$
In fact,
$$\align&\sum_{j=0}^{2^n-1}(-1)^j(2j+1)^k=\sum_{i=0}^{2^n-1}(-1)^{2^n-1-i}\l(2(2^n-1-i)+1\r)^k
\\\eq&-\sum_{i=0}^{2^n-1}(-1)^i\l((2^{n+1}-2i-1)^2\r)^{k/2}
\eq-\sum_{i=0}^{2^n-1}(-1)^i(2i+1)^k\ (\mo\ 2^{n+2}).
\endalign$$
This concludes the proof. \qed

\medskip
\noindent{\it Proof of Theorem 1.2}. Suppose that $k-l=2^nq$ where
$n,q\in\Z^+$ and $2\nmid q$. We want to show $2^n\|(E_k-E_l)$.

 By elementary number theory, $a^{2^{n}}\eq1\ (\mo\
2^{n+2})$ for any odd integer $a$ (cf. [IR, pp.\,43--44]). This, together with
Theorem 1.1, yields that
$$\aligned&\f{3^{k+1}+1}4E_k\eq
\f{3^k}2\sum_{j=0}^{2^{n+1}-1}(-1)^{j-1}(2j+1)^k\l\lfloor\f{3j+1}{2^{n+1}}\r\rfloor
\\\eq&\f{3^l}2\sum_{j=0}^{2^{n+1}-1}(-1)^{j-1}(2j+1)^l\l\lfloor\f{3j+1}{2^{n+1}}\r\rfloor
\eq\f{3^{l+1}+1}4E_l\ (\mo\ 2^{n+1}).
\endaligned$$
As $3^k\eq3^l\ (\mo\ 2^{n+2})$, the odd integers
$(3^{k+1}+1)/4$ and $(3^{l+1}+1)/4$ are congruent modulo $2^n$.
Therefore $E_k\eq E_l\ (\mo\ 2^n)$.

By Theorem 1.1,
$$\f{3^{k+1}+1}4E_k\eq\f{3^k}2\sum_{j=0}^1(-1)^{j-1}
(2j+1)^k\l\lfloor\f{3j+1}2\r\rfloor=3^{2k}\eq1\ (\mo\ 2).$$
So $E_k$ is odd. Similarly, $E_l\eq1\ (\mo\ 2)$. In light of the
above,
$$\align E_k\eq E_l\ (\mo\ 2^{n+1})\iff&\f{3^{k+1}+1}4\eq\f{3^{l+1}+1}4\ (\mo\ 2^{n+1})
\\\iff& 3^{2^nq}\eq1\ (\mo\ 2^{n+3}).
\endalign$$
It is well known that the order of $5$
mod $2^{n+3}$ is $2^{n+1}$ and that $3\eq(-1)^a5^b\ (\mo\
2^{n+3})$ for some $a\in\{0,1\}$ and
$b\in\{0,1,\ldots,2^{n+1}-1\}$ (cf. [IR, pp.\,43--44]). Thus
$$3^{2^nq}\eq1\ (\mo\ 2^{n+3})\iff 5^{2^nbq}\eq1\ (\mo\ 2^{n+3})\iff 2\mid bq\iff 2\mid b.$$
If $2\mid b$, then $(-1)^a3\eq 5^b\eq1\ (\mo\ 8)$ which is
impossible. So $2\nmid b$ and hence $2^{n+1}\nmid(E_k-E_l)$. We are done. \qed

\bigskip
\Ack. The author thanks Prof. Wagstaff for his information on the history of Theorem 1.2
which was rediscovered by the author.

\enddocument